\documentclass [11pt]{article}
\usepackage{amsfonts,amsmath}

\newtheorem{theorem}{Theorem}
\newtheorem{definition}{Definition}
\def\Z{{\mathbb Z}}
\def\R{{\mathbb R}}
\def\C{{\mathbb C}}

\begin{document}

\title{On a high-dimensional generalization of Seifert fibrations}
\author{I.A. Taimanov
\thanks{Sobolev Institute of Mathematics, 630090 Novosibirsk, Russia, and Department of Mechanics and Mathematics, Novosibirsk State University,
630090 novosibirsk, Russia; e-mail: taimanov@math.nsc.ru. \newline The work is supported by RFBR (grant N 12-01-00124-a) and
by the grant of the Goverment of Russian federation for state support of scientific researches (contract N 14.В25.31.0029).}}
\date{}

\maketitle

\hfill{To N.P. Dolbilin on his 70th birthday}

\medskip

%\begin{abstract}
%The notion of generalized Seifert fibration is introduced, it is shown that the projections of certain
%Eschenburg $7$-manifolds $W^7_{\bar{n}}$ onto $\C P^2$ define such fibrations, and for them the characteristic classes corresponding to the generators of
%$H^2(B(U(2)/\Z_{2n});\Z)$ are defined.
%\end{abstract}

\medskip

In this article we consider high-dimensional generalizations of Seifert fibrations which naturally appear in the study of topology of positively curved manifolds
\cite{1} and, in particular, of  certain Eshchenburg $7$-manifolds \cite{Esch}.

We give a rigorous definition of such generalized bundles and discuss the problem of constructing their characteristic classes.

Let $E, M$ and $F$ be smooth manifolds and a Lie group $G$ act on $F$ by homeomorphisms.
Let $G$ contain a normal subgroup which is isomorphic to
$U(1) = \{e^{i\varphi}, \varphi \in \R\}$:
$$
\rho: U(1) \to G.
$$
In the sequel we mean by $G/\Z_k$ the quotient group $G/\rho(\Z_k)$ where $\Z_k \subset U(1)$.

Let $\Sigma$ be a submanifold in $M$ of even codimension and $p,q,r$ be such integers that
$q$ and $r$ are relatively prime.

\begin{definition}
A triple $(E,M,\pi)$, where
$$
\pi: E \to M
$$
is a continuous mapping from $E$ onto $M$, form a generalized Seifert fibration (or, for brevity, $S$-fibration)
with the singular locus $\Sigma$ if

1) the triple $(E \setminus \pi^{-1}(\Sigma), M \setminus \Sigma, \pi)$ is the bundle over $M \setminus \Sigma$ with the fibre
$F/\Z_p$ and the structure group $G/\Z_p$ where the action of $G/\Z_p$ on $F/\Z_p$ is induced by the action of $G$ on $F$;

2) every point $x \in \Sigma$ has a neighborhood $V$ with coordinates $x^1,\dots,x^{2l}$, $x^{2l+1}$, $\dots,x^n$ such that

а) $\Sigma$ is defined by the equations
$$
x^1=\dots=x^{2l}=0;
$$

б) on $W = V \times U(1) \times F$ there is a free $U(1)$--action:
$$
(X,Y,u,v) \stackrel{e^{i\varphi} \in U(1)}{\longrightarrow} (e^{ip\varphi} \cdot X, Y, e^{iq\varphi} \cdot u, e^{ir\varphi} \cdot v),
$$
where
$$
X = (x^1,\dots,x^{2l}) = (z^1=x^1+ix^2,\dots,z^l = x^{2l-1}+ix^{2l}),
$$
$$
e^{i\varphi} \cdot X = (e^{i\varphi}z^1,\dots,e^{i\varphi}z^l),
\ \ \ Y = (x^{2l+1},\dots,x^n),
$$
$$
u \in U(1), \ \ e^{i\varphi}\cdot u = e^{i\varphi}u, \ \ \ v \in F, \ \ e^{i\varphi} \cdot v = \rho(e^{i\varphi})(v) \in F;
$$

в) $\pi^{-1}(V)$ is homeomorphic to the quotient-space of $W$ with respect to the $U(1)$-action and
$$
\pi([(X,Y,u,v)]) = (X,Y) \in V;
$$

3) the triple $(\pi^{-1}(\Sigma),\Sigma,\pi)$ is a fibre bundle over $\Sigma$ with the fibre $F/\Z_q$ and the structure group $G/\Z_q$
where the action of $G/\Z_q$ on $F/\Z_q$ is induced by the action of $G$ on $F$.
\end{definition}

{\sc Remark.} To avoid excessive bulkiness we give a definition in the simplest situation the branch submanifold
$\Sigma$ is either connected, either for all its components the dimension and the constants $p,q$ are the same. This deifnition is naturally generalized 
to the case when $\Sigma$ consists of few components:
$$
\Sigma = \Sigma_1 \cup \dots \cup \Sigma_k,
$$
and, if $F$ is not homeomorphic to the circle $S^1$, then the corresponding parameters $(p_1,q_1,r_1), \dots, (p_k,q_k,r_k)$
meet the condition
$$
p_1 = \dots = p_k.
$$

{\sc Examples.}

1) {\sc Seifert fibrations}. Let $M$ be a $2$-manifold, $\Sigma = Q_1 \cup \dots \cup Q_k$ be a union of points
$Q_1,\dots,Q_k \in M$, $F = S^1 = U(1)$,
$G = U(1)$ act on fibers by left translations, $\rho$ be an isomorphism, $p_j$ are $q_j$ relatively prime and $r_j=1$ for all
$j=1,\dots,k$.
Then $E$ is a classical Seifert fibration \cite{S,Scott}.

2) {\sc Eschenburg $7$-manifolds.}. Let us consider the matrix group $SU(2)$ embedded into $SU(3)$:
$$
A \to \left(\begin{array}{cc} A & 0 \\ 0 & 1 \end{array}\right), \ \ \ A \in SU(2),
$$
and the subgroup $T_{k,l}$, in $SU(3)$, isomorphic to $U(1)$:
$$
T_{k,l} = \left\{ \left(\begin{array}{ccc} e^{ik\varphi} & 0 & 0 \\ 0 & e^{il\varphi} & 0 \\ 0 & 0 & e^{-i(k+l)\varphi} \end{array}\right),
\varphi \in \R \right\},
$$
where $k,l \in \Z$, and the greatest common divisor of $k,l$ and $k+l$ is equal to $1$.
The manifold
$$
N_{k,l} = SU(3)/T_{k.l}
$$
is called the Aloff--Wallach manifold and admits a positively curved metric.

The quotient space $SU(3)/SU(2)$ is diffeomorphic to the $5$-sphere $S^5$ and he projection 
$$
\pi: SU(3) \to S^5 = SU(3)/SU(2)
$$
defines the bundle over $S^5$ with the fibre $S^3$.

Remark that since
$\pi_4(U(2)) = \pi_4(S^3) = \Z_2$ there are two nonisomorphic bundles over $S^5$ with the fibre diffeomorphic to $S^3$ and the structural group
$U(2)$:
the bundle $SU(3) \stackrel{S^3=SU(2)}{\longrightarrow} S^5$ mentioned above and the trivial bundle
$S^5 \times S^3  \stackrel{S^3}{\longrightarrow} S^5$.

For $k=l=1$ the subgroups $SU(2)$ and $T_{1,1}$ commute and hence $T_{1,1}$ acts on $S^5$, the corresponding orbit space is
$\C P^2$ and we obtain the bundle
\begin{equation}
\label{aloff}
N_{1,1} \stackrel{\R P^3}{\longrightarrow} \C P^2,
\end{equation}
(see \cite{1}).

In analogy with (\ref{aloff}) we propose in \cite{1} to consider the fibering mappings of the Eschenburg spaces
$W^7_{\bar{n}}$ onto $\C P^2$. The Eschenburg space $W^7_{\bar{n}}$
пis obtained by the factorization of $SU(3)$ under the free two-side action of $U(1)$:
$$
W^7_{\bar{n}} = T_{2n,2n,2} \setminus SU(3) / T_{-(1+2n),-(1+2n),0},
$$
where
$$
T_{k,l,m} = \left\{ \left(\begin{array}{ccc} e^{ik\varphi} & 0 & 0 \\ 0 & e^{il\varphi} & 0 \\ 0 & 0 & e^{im\varphi} \end{array}\right),
\varphi \in \R \right\}.
$$
These spaces are particular cases of the Eschenburg which are biquotients of 
$SU(3)$ under the action of $U(1)$ and admit metrics of positive sectional curvature \cite{Esch}. By their construction, they are
nonhomogeneous generalizations of the Aloff--Wallach spaces.

In this case the right-side action of $SU(2)$ on $SU(3)$ commutes with the $T_{-(1+2n),-(1+2n),0}$-action and analogously to (\ref{aloff})
we may construct the projection of $W^7_{\bar{n}}$ onto $\C P^2$.
The space $\C P^2$ is the quotient space of $S^5$ with respect to the $U(1)$-action. In final remarks to 
\cite{1}
\footnote{This work was devoted to totally geodesic embeddings of positively curved $7$-spaces into $13$-dimensional spaces of positive curvature 
(the choice of dimensions is due to known examples) and their relation to the pinching constants of metrics. In the article we obtained the very first results 
on this problem which later was studied in \cite{P,DE,K}.}
we briefly defined this mapping and guessed that it is a fibre bundle. Ziller noted that the preimages of points from
$\C P^1 \subset \C P^2$ and from $\C P^2 \setminus \C P^1$ are different and homeomorphic to $\R P^3/\Z_n$ and $\R P^3$ respectively.
Similar examples were discussed later in \cite{FZ}.

The study of such situations needed a rigorous definition of appearing generalized bundles which would give the procedure of constructing them and correctly pose the problem of their topological classification. The main goal of this article is the definition of this class of fibrations as 
$S$-fibrations.

\begin{theorem}
The projection
$$
\pi: W^7_{\bar{n}} \to \C P^2
$$
defines an $S$-fibration (a generalized Seifert fibration) for which
$$
\Sigma = \C P^1, \ \ \ G = U(2), \ \ \ \ F = S^3 (=SU(2)),
$$
$G=U(2)$ acts on $F$ as on the unit sphere in $\C^2$, the action of $U(1)$ on $F=S^3$ is standard:
$$
(z_1,z_2) \in S^3 \stackrel{e^{i\varphi}}{\longrightarrow} (e^{i\varphi}z_1, e^{i\varphi}z_2), \ \ \ (z_1,z_2) \in \C^2, |z_1|^2 + |z_2|^2 =1,
$$
and
$$
p = 2, \ \ \ q = 2n, \ \ \ \ r= 1-2n,
$$
\end{theorem}

{\sc Proof.} The mapping $\pi: W^7_{\bar{n}} \to \C P^2$ is obtained from the projection
$SU(3) \to S^5$ of the form
$$
\tilde{\pi}(A) =  A \cdot \left(\begin{array}{c} 0 \\ 0 \\ 1 \end{array}\right) = \left(\begin{array}{c} a_{13} \\ a_{23} \\
a_{33} \end{array}\right) = \left(\begin{array}{c} z_1 \\ z_2 \\ z_3 \end{array}\right), \ \ \ |z_1|^2+|z_2|^2+|z_3|^2=1,
$$
by the factorization with respect to the fiber-wise action of $U(1)$ which on the base of the projection acts as follows:
$$
(z_1,z_2,z_3) \to (u^{2n} z_1, u^{2n} z_2, u^2 z_3), \ \ \ u= e^{i\varphi} \in U(1).
$$
This action of $U(1)$ on $S^5$ is not free:

a) it is trivial for  $u^2=1$ (however therewith $T_{-(1+2n),-(1+2n),0}$ acts of the fibre $SU(2)$ of the bundle
$SU(3) \to S^5$ as a reflection and hence for $z_3 \neq 0$ the generic fibre is diffeomorphic to $SU(2)/\{\pm 1\} = \R P^3$;

b) if $z_3=0$ and $u^4 = 1$ then the action on $S^3 = \{z_3 =0\}$ is trivial however on the fibre of the bundle
$SU(3) \to S^5$ we obtain an additional $\Z_2$-action.

The case b) was considered in detail in the thesis of our student N.E.~Russ\-kikh whose results were included into
\cite{Rus}. The fibre over $(z_1,z_2,0) \in S^3 \subset S^5$ of the bundle $SU(3) \to S^3$
consists of matrices of the form
$$
C_{a,b} = \left(\begin{array}{ccc} 0 & \bar{z}_2 & z_1 \\ 0 & -\bar{z}_1 & z_2 \\
1 & 0 & 0 \end{array}\right)\left(\begin{array}{ccc} a & b & 0 \\ -\bar{b} & \bar{a} & 0 \\ 0 & 0 & 1 \end{array}\right),
\ \ |a|^2+|b|^2=1
$$
and $U(1)$ acts on this fibre as follows:
$$
C_{a,b} \longrightarrow \left(\begin{array}{ccc} u^{2n} & 0 & 0 \\ 0 & u^{2n} & 0 \\ 0 & 0 & u^2 \end{array}\right)
C_{a,b} \left(\begin{array}{ccc} u^{-(1+2n)} & 0 & 0 \\ 0 & u^{-(1+2n)} & 0 \\ 0 & 0 & 1 \end{array}\right)
=
$$
$$
=
\left(\begin{array}{ccc} 0 & u^{-2n}\bar{z}_2 & u^{2n} z_1 \\ 0 & -u^{-2n}\bar{z}_1 & u^{2n} z_2 \\
1 & 0 & 0 \end{array}\right)\left(\begin{array}{ccc} u^{1-2n} a & u^{1-2n} b & 0 \\ -u^{2n-1} \bar{b} & u^{2n-1} \bar{a} & 0 \\ 0 & 0 & 1 \end{array}\right).
$$
This implies that the mapping $\pi: W^7_{\bar{n}} \to \C P^2$ defines an $S$-fibration for the parameters
$p,q,r$ mentioned in Theorem.

It is natural to consider the problem of defining the characteristic classes of $S$-fibrations.

Clearly they include the characteristic classes of the bundle
\begin{equation}
\label{locus}
\pi^{-1} (\Sigma) \stackrel{F/\Z_q}{\longrightarrow} \Sigma
\end{equation}
over the singular locus.

Let us recall the well-known facts:

1) principal bundles over $X$ with the structure group $G$ are classified by homotopy classes of continuous mappings from
$X$ to $BG$, the classified space of $G$.
The set of such homotopy classes is denoted by $[X,BG]$.
$BG$ is uniquely defined (up to homotopy equivalence) as the base of the principal $G$-fibration 
$E \to BG$ with $E$ contractible;

2) given a continuous homomorphism $f: H \to G$,  there is the induced mapping
$$
Bf: BH \to BG,
$$
with natural functorial properties;

3) the characteristic classes of the principal $G$-bundle over $X$, correspon\-ding to a homotopy class $[\alpha] \in [X,BG]$,
are the images of the distinguished elements $\nu \in H^\ast(BG;\Lambda)$ under the induced mapping
$$
\alpha^\ast: H^\ast(BG;\Lambda) \to H^\ast(X;\Lambda)
\ \ \ : \ \ \
\alpha^\ast(\nu) \in H^\ast(X;\Lambda).
$$

Therefore for defining the characteristic classes of the bundles (\ref{locus}) we have to compute the cohomology
of the spaces $B(U(2)/\Z_q)$. We restrict ourselves to the low-dimensional cohomology groups.

\begin{theorem}
$$
H^1(B(U(2)/\Z_q);\Z)= 0, \ \ \ \ H^1(B(U(2)/\Z_q);\Z_2)= 0,
$$
$$
H^2(B(U(2)/\Z_q);\Z)= \Z,
$$
$$
H^2(B(U(2)/\Z_q);\Z_2)=
\begin{cases} \Z_2 & \mbox{for $q$ odd}; \\
\Z_2 + \Z_2 & \mbox{for $q$ even}.
\end{cases}
$$
 \end{theorem}

{\sc Proof.} The space $U(2)/\Z_q$ is connected.
Let us consider the covering
$$
\pi: U(2) \stackrel{\Z_q}{\to} U(2)/\Z_q
$$
and the corresponding exact homotopy sequence:
$$
\dots \to \pi_i(\Z_q) \to \pi_i(U(2)) \to \pi_i(U(2)/\Z_q) \to \dots.
$$
It implies that
$$
\pi_i(U(2)) = \pi_i(U(2)/\Z_q) , \ \ \ i \geq 2,
$$
and there is an exact sequence
$$
0 \to \pi_1(U(2))=\Z \stackrel{\pi_\ast}{\longrightarrow} \pi_1(U(2)/\Z_q) \stackrel{j_\ast}{\longrightarrow} \Z_q \to 0.
$$

When considering the fundamental groups of Lie groups we shall always consider loops starting and ending at the unit.

It is known that $\pi_1(U(2)) = \Z$ and this group is generated by the homotopy class $[\eta]$ of the loop
$$
\eta(\varphi) = \left(\begin{array}{cc} e^{i\varphi} & 0 \\ 0 & 1 \end{array}\right).
$$
Let us denote by $[\xi]$ the homotopy class of the loop
\begin{equation}
\label{xi}
\xi(\varphi) = \left(\begin{array}{cc} e^{\frac{i\varphi}{q}} & 0 \\ 0 & e^{\frac{i\varphi}{q}} \end{array}\right).
\end{equation}
$j_\ast([\xi])$ generates the subgroup $\Z_q$ and $\pi_1(U(2)/\Z_q)$ is generated by
$[\xi]$ and $[\eta]:= \pi_\ast([\eta])$ meeting the relation
$$
q [\xi] = 2 [\eta].
$$
Let us consider separately two cases:

1) $q = 2s+1$. Then
$$
u = (2s+1)[\xi] - 2[\eta], \ \ \ \ v= -s[\xi]+[\eta]
$$
generate the lattice
$\Z [\xi] + \Z[\eta]$ and
$$
\pi_1(U(2)/\Z_q) = \{\Z u +\Z v \}/\Z u = \Z;
$$

2) $q = 2s$. Then
$$
u = s[\xi] - [\eta], \ \ \ \ v= (1-s)[\xi]+[\eta]
$$
generate the lattice
$\Z [\xi] + \Z[\eta]$ and
$$
\pi_1(U(2)/\Z_q) = \{\Z u +\Z v \}/2\Z u = \Z + \Z_2.
$$

Since $\pi_1(U(2)/\Z_q)$ is commutative, we have
$$
H_1(U(2)/\Z_q) = \pi_1(U(2)/\Z_q).
$$

It immediately follows from the spectral sequence of the bundle 
$$
U(2)/\Z_q \to E \sim \mathrm{pt} \to B(U(2)/\Z_q)
$$
with the contracted space $E$ that
$$
H_1(B(U(2)/\Z_q))= 0, \ \ \
H_2(B(U(2)/\Z_q)) =
\begin{cases} \Z & \mbox{for $q$ odd};\\
\Z + \Z_2 & \mbox{for $q$ even}.
\end{cases}
$$
However in these dimensions it is enough to use the exact homotopy sequence of this bundle.

Together with the universal coefficients formula that implies Theorem.
Theorem is proved.

Let us introduce now the characteristic classes:

1) if $q=2s+1$, then $H^2(BU(2);\Z)$ is generated by the cohomology class adjoint to $v$.
Let us denote it by
$$
\kappa \in H^2(B(U(2)/\Z_q);\Z),
$$
as well as the corresponding characteristic class;

2) if $q=2s$, then $H^2(BU(2);\Z)$ is again generated by the cohomology class $\kappa$ adjoint to $v$.
However $H^2(B(U(2)/\Z_q);\Z_2)$ is generated by two classes
$$
\kappa_u, \kappa_v \in H^2(B(U(2)/\Z_q);\Z_2),
$$
adjoint to $u$ and $v$ correspondingly.

It is clear that the class $\kappa$ similarly to the first Chern class which it generalizes admits a simple geometrical interpretation and we have

\begin{theorem}
\label{class}
Let the principal bundle $E$ over $\C P^1$ with the structural group $U(2)/\Z_q$  is obtained by gluing two trivial bundles
over two-dimensional discs via the homeomorphism
$\psi$ of the trivial bundles on the boundaries and this homeomorphism is defined by the mapping
$S^1 \to U(2)/\Z_q$ which realizes the element
$a[v]+b[u] \in \pi_1(U(2)/\Z_q)$. Then
$$
\kappa(E) = a[\C P^1].
$$
In particular, if $\psi$ is defined by $\xi^m$, where $\xi$ has the form (\ref{xi}), then
$$
\kappa(E) = \begin{cases} 2m\,[\C P^1] & \mbox{for $q$ odd}; \\
m\,[\C P^1] & \mbox{for $q$ even}.
\end{cases}
$$
\end{theorem}

This theorem implies the following. Let us consider the cohomology groups of $B(U(2)/\Z_q)$ with coefficients in $\Z\left[\frac{1}{N}\right]$, the commutative group generated over 
$\Z$ by $\frac{1}{N}$ where
$$
N = \begin{cases} q & \mbox{for $q$ odd} \\
\frac{q}{2} & \mbox{for $q$ even}.
\end{cases}
$$
The natural homomorphism of division by $N$
$$
\Z \stackrel{\frac{1}{N}}{\longrightarrow} \Z\left[\frac{1}{N}\right]
$$
induces the isomorphism
$$
H^2(B(U(2)/\Z_q);\Z) \to H^2(B(U(2)/\Z_q); \Z\left[\frac{1}{N}\right]),
$$
therewith $\kappa$ goes into a generator of $H^2(B(U(2)/\Z_q; \Z\left[\frac{1}{N}\right])$ over $\Z$ for which we save the same notation $\kappa$.
The subgroup
$$
H^2(B(U(2)/\Z_q);\Z) \subset H^2(B(U(2)/\Z_q); \Z\left[\frac{1}{N}\right])
$$
is isomorphic to $H^2(BU(2);\Z)$ and is generated by
$$
c_1 = N \kappa,
$$
the first Chern class of bundles whose structural group is lifted to $U(2)$.
Therewith $\kappa = \frac{c_1}{N}$ is naturally treated as the {\it fractional Chern class}.

Let expose the main results of the article \cite{Rus} by N.E. Russkikh:

\begin{itemize}
\item
$$
\kappa = (2n-1)[\C P^1] \in H^2(\C P^1) = \Z
$$
for the $S$-fibrations $W^7_{\bar{n}} \to \C P^2$;

\item
there exists an $S$-fibration with the same data $M,\Sigma,F,G,p,q,r$ as of the fibration $W^7_{\bar{2}} \to \C P^2$ such that it is obtained from a fiber-wise action of
$U(1)$ on the trivial bundle $S^5 \times S^3 \to S^5$ such that it has the same value of the characteristic class $\kappa$ as the $S$-fibration $W^7_{\bar{2}} \to \C P^2$.
\end{itemize}

The last example shows that for recognizing $S$-fibrations it is desirable to define higher characteristic classes.

For simplicity let us restrict ourselves to $S$-fibrations over $\C P^2$ with the singular locus $\Sigma = \C P^1$.
In this case one may define the difference element of $S$-fibrations when they are isomorphic as $S$-fibrations over a tubular neighborhood
$U(\Sigma)$ of the singular locus $\Sigma= \C P^1$. The boundary $\partial\Sigma$
is diffeomorphic to the $3$-sphere $S^3$ which from one side bounds a four-dimensional cell $D^4$. Every principal $U(2)/\Z_p$-bundle over $S^3$ is trivial. Fixing
trivializations of the bundles over $\partial U(\Sigma)$ and $\partial D^4$, the $S$-fibration over $U(\Sigma)$ and the trivial bundle over $D^4$ may be glued into
an $S$-fibartion over $\C P^2$ via a homeomorphism of trivial bundles over the boundaries. Up to homotopy such a homeomorphism is defined by an element
$$
\alpha \in \pi_3(U(2)/\Z_q) = \Z.
$$
If gluings correspond to the elements $\alpha_1$ and $\alpha_2$, then we have the difference element
\begin{equation}
\label{dif}
\alpha_1 - \alpha_2 \in H^4(\C P^2,\C P^1;\Z).
\end{equation}
For known examples the difference elements are not computed until recently.

The following questions are staying open:

\begin{enumerate}
\item
how to define the analog of the second Chern class for $S$-fibrations over $(\C P^2, \C P^1)$?

\item
how to define the Euler class for general $S$-fibrations (for Seifert fib\-ra\-tions such a definition exists \cite{Scott})?

\item
up to which preciseness the restriction of an $S$-fibration over a tubular neighborhood of the singular locus
$\Sigma$ is defined by its restriction over $\Sigma$?
\end{enumerate}

{\sc Remark.} The referee draw our attention to the most general definition of Seifert bundles given in 
\cite{H} (see also \cite{O} where it is exposed in \S 5). It is based on general constructions of 
the theory of $G$-spaces and there is the following theorem:

{\sl Let a locally compact topological group $G$ acts on a locally compact space $X$
such that every mapping $g: X \to X$ is proper and all isotropy subgroups are finite.
Then $(X,\pi,X/G)$, where $\pi: X \to X/G$ is the projection, is a principal Seifert bundle with the structural group
$G$.}

\noindent
Indeed the manifolds from Theorem 1 satisfy this definition but for another structural group $G=SU(2)$.
In \cite{H} all concrete examples which are discussed are classical Seifert $3$-manifolds and the problem of defining the characteristic classes is not considered.
For the class of manifolds considered by us the structural group is defined differently and that allows to introduce the charac\-te\-ristic class
$\kappa$ and the difference element (\ref{dif}). Apparently that has to lead to defining higher characteristic classes.

The author thanks N.E. Russkikh for helpful discussions and the referee for helpful remarks.

\end{document}